\documentclass[11pt,twoside,letterpaper]{article}
\author{Christine M. Escher, S.K.Ultman}
\date{6 December 2008}
\title{Cohomology rings of certain seven dimensional manifolds}

\usepackage{textcomp} 
\usepackage{amsmath,amssymb}  
\usepackage{mathtools}
\usepackage{eufrak} 
\usepackage[mathcal]{euscript} 
\usepackage{mathrsfs} 
\usepackage{bbm}  
\usepackage[all]{xy} 
\usepackage{amscd} 

\newtheorem{thm}{Theorem}[section]
\newtheorem{prop}[thm]{Proposition} 
\newtheorem{cor}[thm]{Corollary} 
\newtheorem{lemma}[thm]{Lemma}

\def\Z{\mathbb{Z}}                       
\def\R{\mathbb{R}}

\def\H{\mathbb{H}}
\def\1{\mathbbm{1}}
\def\qed {\hspace{\stretch{1}} $\square$\\}  
\def\T{\rule{0pt}{3ex}}                 
\def\B{\rule[-1.5ex]{0pt}{0pt}}        

\begin{document}

\maketitle

\begin{abstract}
We calculate the cohomology rings of a collection of seven dimensional manifolds supporting an $S^3 \times S^3$-action with one dimensional orbit space.  These manifolds are of interest to differential geometers studying non-negative and positive sectional curvature.  From this collection, we identify several families of manifolds for which there exist well-known topological invariants distinguishing homeomorphism and diffeomorphism types.
\end{abstract}

\section{Introduction and background.} \label{section1}

In the study of Riemannian manifolds of positive sectional curvature, differential geometers face the following challenge:  although there are few known examples, only a small number of obstructions to positive curvature have been found.  Recently, a new example of positive curvature was discovered in an infinite family of seven dimensional manifolds (\cite{Dearricott}, \cite{GVZ}).  Related to this family of manifolds are three other infinite families of seven dimensional manifolds; in this article, we determine the cohomology rings of these three families.  The cohomology rings are then used to identify manifolds which possess an invariant of diffeomorphism type, which is computable using purely topological means (as in \cite{KS}).

The manifolds under consideration are of cohomogeneity one, where a smooth manifold $M$ is of cohomogeneity one if it is closed, connected, and supports a smooth action by a compact Lie group $G$ with one dimensional orbit space $M/G$.  Cohomogeneity one actions on simply connected manifolds in dimensions five through seven were classified in \cite{Hoel}.  This article is concerned with four of the seven dimensional families resulting from this classification, listed here in Table \ref{table1}.  By \cite{GZ} Theorem E, it is known that members of the families $L_{(p_-,q_-)(p_+,q_+)}$, $M_{(p_-,q_-)(p_+,q_+)}$ and $N_{(p_-,q_-)(p_+,q_+)}$ admit non-negative curvature; however, it is unknown whether this is true in general for the family $O_{(p,q:m)}$.  The families $M_{(p_-,q_-)(p_+,q_+)}$ and $N_{(p_-,q_-),(p_+,q_+)}$ are of particular interest due to the fact that some members may admit positive sectional curvature (see \cite{GWZ}); the new example of positive curvature mentioned above is in fact a member of the family $M_{(p_-,q_-)(p_+,q_+)}$.  The cohomology groups for members of the families $M_{(p_-,q_-)(p_+,q_+)}$ and $N_{(p_-,q_-)(p_+,q_+)}$ were calculated in \cite{GWZ}.  Although the cohomology rings of members of $M_{(p_-,q_-)(p_+,q_+)}$ are determined by the cohomology groups, the rings of members of $N_{(p_-,q_-)(p_+,q_+)}$ are not obvious.  They are described here for the first time, together with the cohomology rings of members of $L_{(p_-,q_-)(p_+,q_+)}$ and $O_{(p,q:m)}$ having non-trivial finite cyclic fourth cohomology group.   

\begin{table}[t]
\begin{tabular}{@{}c@{\;}|@{\;}c@{}}
\hline \hline
\T \B Family & Isotropy groups $H \subseteq K_-,\,K_+$\\
\B & and restrictions on the parameters\\
\hline
\T \B $L_{(p_-,q_-),(p_+,q_+)}$ & $\langle (i,i) \rangle \subseteq \{(e^{ip_-\theta},e^{iq_-\theta})\},\,\{(e^{jp_+\theta},e^{jq_+\theta})\} \cdot H$\\
\B & $p_-, q_- \equiv 1  \textrm{ mod }4$ \\
\hline
\T \B $M_{(p_-,q_-),(p_+,q_+)}$ & $\Delta Q \subseteq \{(e^{ip_-\theta},e^{iq_-\theta})\} \cdot H,\,\{(e^{jp_+\theta},e^{jq_+\theta})\} \cdot H$\\ 
\B & $\Delta Q$ the diagonal embedding of $\langle 1,i,j,k \rangle$;\\
\B & $p_\pm, q_\pm \equiv 1  \textrm{ mod }4$ \\
\hline
\T \B $N_{(p_-,q_-),(p_+,q_+)}$ & $\langle (h_1,h_2),(1,-1) \rangle \subseteq \{(e^{ip_-\theta},e^{iq_-\theta})\} \cdot H,\,\{(e^{jp_+\theta},e^{jq_+\theta})\} \cdot H$\\ 
\B & $h_1,h_2 \in \{i,-i\}$ with signs chosen so that\\
\B & $(h_1,h_2)$ lies in $\{(e^{ip_-\theta},e^{iq_-\theta})\}$;\\
\B & $p_-$ and $q_\pm$ odd, $p_+$ even\\
\hline
\T \B $O_{(p,q:m)}$ & $\Z_m \subseteq \{(e^{ip\theta},e^{iq\theta})\},\,\Delta S^3 \cdot H$\\
\B & either $m=1$ (with no restrictions on $p$ or $q$) \\
\B & or $m=2$ and $p$ is even\\
\hline \hline
\end{tabular}
\caption{Isotropy groups description of simply connected seven dimensional primitive cohomogeneity one manifolds.} \label{table1}
\end{table}

Our manifolds belong to the class of cohomogeneity one manifolds whose orbit spaces are diffeomorphic to a closed interval.  For this class of manifolds, the additional topological structure induced by the group action has a particularly nice description.  Specifically, when $M/G$ is diffeomorphic to a closed interval, there are (up to conjugation in $G$) two non-principal isotropy groups $K_\pm$ and a principal isotropy group $H$ of the $G$-action such that the non-principal orbits over the endpoints of the interval are diffeomorphic to $G/K_\pm$, while a principal orbit over an interior point is diffeomorphic to $G/H$.  In fact, these groups are sufficient to describe such a manifold.  A manifold $M$ supports a cohomogeneity one action by a compact Lie group $G$ if and only if there are closed subgroups and inclusions $H \subseteq K_-, K_+ \subseteq G$ such that $K_\pm/H$ are diffeomorphic to spheres $S^{t_\pm - 1}$ and $M$ is diffeomorphic to the union of the total spaces $D(G/K_\pm) = G \times_{K_\pm} D^{t_\pm}$ of disk bundles over $G/K_\pm$.  These bundles are glued together along their common boundary $G/H$, which is diffeomorphic to the total space $G \times_{K_\pm} (K_\pm/H) \approx G \times_{K_\pm} S^{t_\pm - 1}$ of the boundary sphere bundles of the disk bundles.  The principal orbits of this action are diffeomorphic to $G/H$.  The non-principal orbits are diffeomorphic to the bases of the disk bundles $D(G/K_\pm)$, and are identified with the zero sections in their respective total spaces $D(G/K_\pm)$.  For further details, see \cite{GZ}.

This description of $M$ is an example of a more general topological construction, which we will call a double disk bundle.  This is a quotient space $X = D(B_-) \cup_\varphi D(B_+)$, where the two spaces $D(B_\pm)$ are the total spaces of disk bundles over paracompact bases $B_\pm$ and the attaching map $\varphi$ is a homeomorphism of the boundaries $\partial D(B_\pm)$.  In the case of a cohomogeneity one manifold, the spaces $B_\pm$ are the non-principal orbits $G/K_\pm$, the boundaries $\partial D(G/K_-) = \partial D(G/K_+) = G/H$, and the attaching map $\varphi$ is the identity map.  So for $M$ cohomogeneity one, $M = D(G/K_-) \cup_{id} D(G/K_+)$.  This double disk description allows us to define long exact sequences which prove useful in determining the cohomology of the manifolds.

All manifolds under consideration admit a cohomogeneity one action by $G = S^3 \times S^3$.  We regard $S^3$ as the group of unit quaternions.  Given a circle group $\{(e^{ip\theta}, e^{iq\theta})\} \subset S^3 \times S^3$, it is always assumed that the parameters $p$ and $q$ are relatively prime integers.  Table \ref{table1} lists the remaining data necessary to describe the manifolds; namely, the isotropy groups together with restrictions on the parameters that guarantee an embedding of the principal isotropy group $H$ in the non-principal isotropy groups $K_\pm$.  For example, in the description of the family $L_{(p_-,q_-),(p_+,q_+)}$ the principal isotropy group $H = \langle(i,i)\rangle$ is the cyclic group of order four generated by the diagonal embedding of the unit quaternion $i$ (the notation $\langle q_1,\ldots,q_n\rangle$ denoting the subgroup generated by the elements $q_1,\ldots,q_n$), the non-principal isotropy group $K_+ = \{(e^{jp_+\theta}, e^{jq_+\theta})\} \cdot H$ is the group whose elements are products of an element of the circle group with an element of $H$, and the congruence of the parameters $p_-, q_-  \equiv 1  \textrm{ mod }4$ of the non-principal isotropy group $K_- = \{(e^{ip_-\theta}, e^{iq_-\theta})\}$ ensures that $H$ is a subgroup.  In comparing this table with Table I of \cite{Hoel}, one should observe that the original notation of \cite{Hoel} has been modified slightly in order to agree with \cite{GWZ}.

Let $M$ be a smooth, simply connected manifold.  Following \cite{EZ}, we will say $M$ has cohomology type $E_r, \, r \geq 2$, if it has non-trivial integral cohomology groups $H^0(M) = H^2(M) = H^5(M) = H^7(M) \cong \Z$ and $H^4(M) \cong \Z_r$; and if the square of a generator of $H^2(M)$ generates $H^4(M)$.  In \cite{KS}, Kreck and Stolz developed topological invariants that can be used to classify these manifolds up to homeomorphism and diffeomorphism.  Due the existence of these powerful invariants, our main result identifies the manifolds listed in Table \ref{table1} which are of cohomology type $E_r$:

\begin{prop} \label{typeErprop}
A simply connected seven dimensional primitive cohomogeneity one manifold has cohomology type $E_r$ (for some $r \geq 2$) if and only if:
\begin{enumerate}
\item[a)] it is a member of the family $L_{(p_-,q_-),(p_+,q_+)}$ with the parameter $p_+$ odd and $p_+^2q_-^2 - p_-^2q_+^2 \neq 0$, or:
\item[b)] it is a member of the family $N_{(p_-,q_-),(p_+,q_+)}$, or:
\item[c)] it is a member of the family $O_{(p,q:m)}$ with $|p|$ and $|q|$ not both equal to one.
\end{enumerate}
\end{prop}

An interesting class of manifolds having cohomology type $E_r$ are the Eschenburg spaces, which are seven dimensional biquotients of $SU(3)$ under a free two-sided circle action.  They have cohomology type $E_r$ where $r$ is odd, and cohomology rings generated by generators of the second and fifth cohomology groups (see \cite{Kruggel1} and \cite{Eschenburg2}).  All such manifolds are known to admit non-negative sectional curvature (see \cite{Eschenburg1}), and some have been shown to admit positive sectional curvature.  A diffeomorphism classification of these spaces has been carried out in \cite{Kruggel2}, \cite{AMP} and \cite{CEZ} using the Kreck-Stolz invariants.  Thus, the Eschenburg spaces provide a good ``comparison'' family for manifolds of cohomology type $E_r$ with $r$ odd.  In light of this, an interesting corollary of Proposition \ref{typeErprop} is:

\begin{cor} \label{eschencor}
A simply connected seven dimensional primitive cohomogeneity one manifold has the integral cohomology ring of an Eschenburg space if and only if:
\begin{enumerate}
\item[a)] it is any member of the family $N_{(p_-,q_-),(p_+,q_+)}$ with $|p_-^2q_+^2 - p_+^2q_-^2| \neq 1$, or:
\item[b)] it is a member of the family $O_{(p,q:m)}$ and one of the parameters $p$ or $q$ is even.
\end{enumerate}
\end{cor}
It should be noted that the manifolds $O_{(p,p \pm 1:2)}$ are known to be Eschenburg spaces (see \cite{GWZ}).

In Section \ref{section2}, we introduce two long exact cohomology sequences and conditions on a general double disk bundle that can be used to determine certain generators of its cohomology ring (Lemma \ref{genlemma}).  Cohomology groups are calculated in Section \ref{section3}, and we determine what conditions on the parameters are necessary to ensure that the fourth cohomology group is non-trivial and finite cyclic.  We note which manifolds have cohomology groups in agreement with those of a manifold having cohomology type $E_r$, and what additional restrictions on the parameters guarantee that the order $r$ of the fourth cohomology group is odd.  In Section \ref{section4}, cohomology rings are computed, completing the proofs of Proposition \ref{typeErprop} and Corollary \ref{eschencor}.  Unless explicitly stated otherwise, integral cohomology is assumed.  We take the cyclic group $\Z_r$ to be infinite cyclic if $r = 0$ and trivial if $r = 1$.

The second author will incorporate the work in this article into a thesis, to be submitted in 2009 in partial fulfillment of the requirements for the degree of Doctor of Philosophy at Oregon State University.

\section{Cohomology of double disk bundles.} \label{section2}

The cohomology groups of a double disk bundle $X = D(B_-) \cup_\varphi D(B_+)$ can be computed in terms of the disk bundles $D(B_\pm) \rightarrow B_\pm$ using the Mayer-Vietoris sequence or the long exact sequences of the pairs $(X,B_\pm)$.  If (as is the case with a cohomogeneity one manifold) the attaching map $\varphi$ is the identity map, the Mayer-Vietoris sequence can be modified as follows.  Let $\partial D(B)$ be the common boundary of the bundles $D(B_\pm)$.  Note that $\partial D(B) \xrightarrow{ \pi_\pm} B_\pm$ are sphere bundles with projections the restrictions of the projections of the disk bundles, and define the homomorphism $ \pi^* := \pi_-^* - \pi_+^*$ to be the difference of the homomorphisms induced by these restricted projections.  The deformation retractions of $D(B_\pm)$ onto $B_\pm$ induce isomorphisms of the cohomology groups.  Making the appropriate substitutions in the Mayer-Vietoris sequence gives the long exact sequence:
\begin{equation} \label{eq:mv}
\cdots \rightarrow H^k(X) \xrightarrow{\psi} H^k(B_-) \oplus H^k(B_+) \xrightarrow{ \pi^*} H^k(\partial D(B)) \xrightarrow{\delta} H^{k+1}(X) \rightarrow \cdots
\end{equation}
where $\psi$ is the composition of the homomorphisms induced by the inclusions of $D(B_\pm)$ in $X$ with the deformation retractions of $D(B_\pm)$ onto $B_\pm$.  The homomorphism $\delta$ is the boundary homomorphism of the Mayer-Vietoris sequence.

The long exact sequence of the pair $(X,B_+)$ can also be modified, assuming the disk bundle $D^t \hookrightarrow D(B_-) \rightarrow B_-$ is orientable; that is, if the structure group $O_t(\R)$ can be reduced to $SO_t(\R)$.  Bundle orientability guarantees the existence of an integral Thom isomorphism from $H^{k-t}(B_-)$ to $H^k(D(B_-), \partial D(B))$.  The bundle projection followed by the inclusion of $B_-$ in $X$ is homotopic to the inclusion of $D(B_-)$ in $X$; so by \cite{Dold} Corollary 11.20, the composition of the Thom isomorphism with the inverse of the excision isomorphism is an $H^*(X)$-module homomorphism from $H^{k-t}(B_-)$ to $H^k(X, D(B_+))$.  The inclusion of the pair $(X,B_+)$ in the pair $(X,D(B_+))$ induces isomorphisms on the relative cohomology groups, as can be seen by considering the long exact cohomology sequences of the pairs and applying the five lemma.  Hence, there is an $H^*(X)$-module isomorphism from the cohomology group $H^{k-t}(B_-)$ to $H^k(X,B_+)$.  Define $J$ to be the composition of this isomorphism with the homomorphism from $H^k(X,B_+)$ to $H^k(X)$ in the long exact sequence of the pair $(X,B_+)$.  One nice property of $J$ is that it is an $H^*(X)$-module homomorphism; this can be very useful when trying to identify generators of the cohomology ring of $X$.

Making the appropriate substitutions in the sequence of the pair $(X,B_+)$, we have the long exact sequence:
\begin{equation} \label{eq:les}
\cdots \rightarrow H^{k-t}(B_-) \xrightarrow{J} H^k(X) \xrightarrow{i_+^*} H^k(B_+) \xrightarrow{\delta}H^{k-t+1}(B_-) \rightarrow \cdots
\end{equation}
(compare with Sequences $4.1.a$ and $4.1.b$ in \cite{Hebda}).  An analogous sequence, with the roles of $B_+$ and $B_-$ reversed, exists if the bundle $D(B_+) \rightarrow B_+$ is orientable.  Furthermore, such a sequence exists for any bundle, regardless of orientability, if integral coefficients are replaced by $\Z_2$-coefficients.  

If $B$ is a closed orientable submanifold of an orientable manifold $M$, then it is known that the normal disk sub-bundle over $B$ in the tangent bundle of $M$ is an orientable bundle (see, for example, \cite{BT} p.66).  As we shall see, both non-principal orbits for manifolds in the family $O_{(p,q:m)}$ are orientable, while only the orbit $G/K_-$ is orientable for members of $L_{(p_-,q_-),(p_+,q_+)}$ and $N_{(p_-,q_-),(p_+,q_+)}$.  Since the orbits are closed submanifolds of a simply connected manifold, at least one long exact sequence of this type exists for each of these manifolds.  Note that both non-principal orbits of members of $M_{(p_-,q_-),(p_+,q_+)}$ are non-orientable (see \cite{GWZ}).

A direct consequence of Sequences \eqref{eq:mv} and \eqref{eq:les} is the following:
\begin{lemma} \label{cycliclemma}
Let $X = D(B_-) \cup_{id} D(B_+)$ be a double disk bundle where the bundle $D^t \hookrightarrow D(B_-) \rightarrow B_-$ is orientable.  For a fixed integer $ \kappa $, suppose  $H^{ \kappa -t}(B_-)$ is cyclic and both groups $H^ \kappa (B_\pm)$ are trivial.  Suppose $H^{ \kappa -1}(\partial D(B))$ is a free group having the same rank as the free part of $H^{ \kappa -1}(B_-) \oplus H^{ \kappa -1}(B_+)$, and that the determinant of the restriction of the homomorphism $ \pi^*$ to the free part of $H^{ \kappa -1}(B_-) \oplus H^{ \kappa -1}(B_+)$ has absolute value $r \geq 0$.  Then $H^ \kappa (X)$ is the finite cyclic group $\Z_r$.
\end{lemma}
\smallskip
\textit{Proof.}  Let $k = \kappa $ in Sequence \eqref{eq:les} to see that $H^ \kappa (X)$ must be a cyclic group.  Sequence \eqref{eq:mv} can be used to show $H^ \kappa (X)$ is equal to the cokernel of the restriction of $ \pi^*$ by setting $k$ equal to $\kappa -1$ and $\kappa$, and observing that $H^{ \kappa}(B_-) \oplus H^{ \kappa}(B_+)$ is trivial.  The Smith normal form of the restriction of $ \pi^*$ can be used to determine the order of $H^\kappa(X)$.  \qed

Given the computability of the cohomology groups of a double disk bundle with respect to its component disk bundles, it seems natural to wonder to what extent these bundles can be used in determining the cohomology ring structure.  The following lemma is a first step in addressing this question:
\begin{lemma} \label{genlemma}
Let $X = D(B_-) \cup_\varphi D(B_+)$ be a double disk bundle over a connected base, where the disk bundle $D^t \hookrightarrow D(B_-) \rightarrow B_-$ is orientable.  Suppose $H^t(X)$ is infinite cyclic, and $H^t(B_+)$ is finite cyclic of order $n \geq 1$.  Let $i_\pm^*$ be the homomorphisms induced on cohomology by the inclusions of $B_\pm$ in $X$, and suppose $i_+^*:H^t(X) \rightarrow H^t(B_+)$ is a surjection.  Finally, suppose $\kappa$ is a fixed integer, $\kappa > t$, such that the following hold:
\begin{enumerate}
\item \label{condition1}  $H^\kappa(X)$ is a non-trivial cyclic group and $H^\kappa(X) \xrightarrow{i_+^*} H^\kappa(B_+)$ is the zero homomorphism.
\item \label{condition2}  $H^{ \kappa -t}(B_-) \cong \Z \cdot \gamma \oplus T$ where $T$ is torsion and the free part is generated by $\gamma$.  If $H^\kappa(X)$ is finite, the orders of elements of $T$ are relatively prime to the order of $H^\kappa(X)$.
\item \label{condition3}   There exists a class $\alpha$ in $H^{\kappa -t}(X)$ with image $i_-^*(\alpha) =  s \gamma + \beta$ (for $\beta \in T$) such that:  if $H^\kappa(X)$ is free, then $|s| = n$; otherwise, $s$ is relatively prime to the order of $H^\kappa(X)$.
\end{enumerate}
Then the cohomology class $x \smile \alpha$ generates $H^ \kappa (X)$, where $x$ is a generator of $H^t(X)$.
\end{lemma}
\smallskip
\textit{Proof}  Let $\1 _-$ be the unit of the cohomology ring $H^*(B_-)$.  Setting $k = t$ in Sequence \eqref{eq:les}, and assuming the hypotheses regarding $H^t(X)$ and $H^t(B_+)$ hold, one has a short exact sequence:
\begin{displaymath}
0 \rightarrow H^0(B_-) \cong \Z \xrightarrow{J} H^t(X) \cong \Z \xrightarrow{i_+^*} H^t(B_+) \cong \Z_n \rightarrow 0.
\end{displaymath}
From this we conclude that the homomorphism $J$ from $H^0(B_-)$ to $H^t(X)$ is multiplication by $n$, and $J(\1 _-) = \pm nx$.

Now, let $k = \kappa$ in Sequence \eqref{eq:les}.  By Condition \ref{condition1},  the homomorphism from $H^\kappa(X)$ to $H^\kappa(B_+)$ is the zero homomorphism, hence by exactness the homomorphism $J$ from  $H^{ \kappa -t}(B_-) \cong \Z \cdot \gamma \oplus T$ to $H^ \kappa (X)$ is a surjection.  Torsion elements of $H^{\kappa -t}(B_-)$ are in the kernel of $J$ (by Condition \ref{condition2}), so $J(\gamma)$ generates $H^\kappa(X)$.

We now consider separately the case in which $H^ \kappa (X)$ is infinite cyclic, and that in which it is finite cyclic.  First, suppose $H^\kappa(X)$ is infinite cyclic.  Let $i_-^*(\alpha) =  \pm n\gamma + \beta \in H^{ \kappa -t}(B_-)$ where $\beta$ is torsion, as required by Condition \ref{condition3}.  Then:
\begin{displaymath}
\pm nJ(\gamma) = \pm (J(n\gamma) + J(\beta)) = \pm J(n\gamma + \beta) = \pm J(i_-^*(\alpha)) = \pm J(\1 _- \smile i_-^*(\alpha)).
\end{displaymath}
Recall that $J$ is an $H^*(X)$-module homomorphism, so:
\begin{displaymath}
J(\1 _- \smile i_-^*(\alpha)) =  J(\1 _-) \smile \alpha = \pm n(x \smile \alpha).
\end{displaymath}
Since the generator $J(\gamma)$ is non-trivial in $H^ \kappa (X) \cong \Z$ and $n \neq 0$, cancellation implies that $x \smile \alpha = \pm J(\gamma)$; so $x \smile \alpha$ generates $H^ \kappa (X)$.

On the other hand, suppose $H^ \kappa (X)$ is finite cyclic.  Let $i_-^*(\alpha) = s\gamma + \beta$ where $s$ is relatively prime to the order of $H^\kappa(X)$, thus satisfying Condition \ref{condition3}.  A calculation similar to the one carried out in the previous case shows that $sJ(\gamma) = \pm n(x \smile \alpha)$.  The class $J(\gamma)$ generates $H^ \kappa (X)$, and the order of $H^ \kappa (X)$ is relatively prime to $s$, therefore $sJ(\gamma) = \pm n(x \smile \alpha)$ also generates $H^ \kappa (X)$.  But if a multiple of $x \smile \alpha$ generates a finite cyclic group, then $x \smile \alpha$ itself must be a generator.  \qed

A useful tool for determining whether Condition \ref{condition3} of Lemma \ref{genlemma} hold, is the commutative ladder of long exact sequences:
\begin{equation} \label{eq:ladder} \xymatrix{
\cdots  \ar[r]^{j_-^* \quad}  & H^k(X) \ar[r]^{i_-^* \quad} \ar[d]_{i_+^*}  & H^k(B_-) \ar[d] \ar[r]^{\delta_- \qquad}       & H^{k+1}(X, B_-) \ar[d]_{\cong} \ar[r]^{\qquad \qquad j_-^*}      & \cdots \\
\cdots \ar[r]^{j^* \quad}      & H^k(B_+)  \ar[r]^{i^* \quad}           & H^k(\partial D(B)) \ar[r]^{\delta \qquad} & H^{k+1}(D(B_+),\partial D(B)) \ar[r]^{\qquad \qquad j^*}           & \cdots.}
\end{equation}
This diagram arises as follows:  begin with the commutative ladder of the long exact sequences of the pairs $(X,B_-)$ and $(X,D(B_-))$, induced by the inclusion of $(X,B_-)$ in $(X,D(B_-))$.  By the five lemma, all vertical homomorphisms in this ladder are isomorphisms (since the inclusion of $B_-$ in $D(B_-)$ is a homotopy equivalence).  The long exact sequence of $(X,B_-)$ forms the top row of the above diagram.  To construct the remainder of the diagram, we use the commutative ladder of long exact sequences induced by the inclusion $e$ of the pair $(D(B_+), \partial D(B))$ in $(X,D(B_-))$.   Note that the induced homomorphism $e^*$ of the relative cohomology groups is the excision isomorphism.  This gives a commutative ladder between the long exact sequences of the pairs $(X,B_-)$ and $(D(B_+), \partial D(B))$ where the vertical homomorphism between the relative cohomology groups is an isomorphism.  Finally, replace the cohomology groups of $D(B_+)$ with those of $B_+$ (allowable, since the inclusion of $B_+$ in $D(B_+)$ is a homotopy equivalence) to arrive at Diagram \eqref{eq:ladder}.

\section{Proof of Proposition 1.1:  the cohomology groups.} \label{section3}

A manifold $M$ having cohomology type $E_r$ has nontrivial cohomology groups $H^0(M) \cong H^2(M) \cong H^5(M) \cong H^7(M) \cong \Z$ and $H^4(M) \cong \Z_r$ a non-trivial finite cyclic group.  If $M$ is an Eschenburg space, $r$ is odd.  The initial step in the proof of Proposition \ref{typeErprop} will be to calculate the cohomology groups of the families $L_{(p_-,q_-),(p_+,q_+)}$, $M_{(p_-,q_-),(p_+,q_+)}$, $N_{(p_-,q_-),(p_+,q_+)}$ and $O_{(p,q:m)}$.  Observe that all manifolds in question are closed, seven dimensional and simply connected (see \cite{Hoel}); as such, they have infinite cyclic cohomology in dimensions zero and seven and trivial cohomology in dimensions one and six.

\subsection{Members of the family $L_{(p_-,q_-),(p_+,q_+)}$.}

Recall that this family is described by the groups:
\begin{displaymath}
H = \langle (i,i) \rangle \subseteq K_- = \{(e^{ip_-\theta},e^{iq_-\theta})\}, K_+ = \{(e^{jp_+\theta},e^{jq_+\theta})\} \cdot H \subseteq G = S^3 \times S^3
\end{displaymath}
where $p_-,q_-$ and $p_+,q_+$ are pairs of relatively prime integers, and $p_-$ and $q_-$ are both congruent to 1 modulo 4.  This family naturally splits into two subfamilies, depending on whether $p_+$ is even or odd.  The cohomology of the principal orbit $G/H$ and the non-principal orbit $G/K_-$ is the same in both cases.  The principal orbit $G/H = S^3 \times S^3/\langle (i,i) \rangle$ is homeomorphic to the product $S^3 \times (S^3/\langle i \rangle)$ of the 3-sphere with the lens space $S^3/\langle i \rangle \approx L_4(1,1)$, with an explicit homeomorphism given by $[q_1,q_2] \mapsto (q_1{q_2}^{-1}, [q_2])$.  The non-principal orbit $G/K_- = S^3 \times S^3 / \{(e^{ip_-\theta},e^{iq_-\theta})\}$ is always homeomorphic to $S^3 \times S^2$ by \cite{WZ} Proposition 2.3.  The orbit $G/K_+$, however, varies depending on the parity of $p_+$.

\smallskip
\noindent \textbf{Case 1.}  Suppose $p_+$ is odd.  Then the cohomology groups of the non-principal orbits $G/K_+$ were calculated in \cite{GWZ} Lemma 13.3a, where they were shown to be:
\begin{displaymath}
H^k(G/K_+) \cong \left\{ \begin{array}{ll}
\Z            & k = 0,3\\
\Z_2       & k = 2,5\\
0             & \textrm{else}.
\end{array} \right.
\end{displaymath}

Let $L$ be a member of this subfamily.  We know that $H^0(L) \cong H^7(L) \cong \Z$.  The orbit $G/K_- \approx S^3 \times S^2$ is a closed orientable submanifold of codimension $2$, so the normal disk bundle over $G/K_-$ is an orientable bundle with fiber $D^2$.  Setting $t=2$, $\kappa = 4$, $B_\pm = G/K_\pm$ and $\partial D(B) = G/H$, it follows from Lemma \ref{cycliclemma} that $H^4(L) \cong coker(\pi^*) \cong \Z_r$.  Recall that $r$ is (up to sign) the determinant of the homomorphism $\pi^*$ from the rank two free abelian group $H^3(G/K_-) \oplus H^3(G/K_+) \cong \Z \oplus \Z$ to the rank two free abelian group $H^3(G/H) \cong \Z \oplus \Z$.  Apply Sequences \eqref{eq:mv} and \eqref{eq:les} (taking $t=2$) to find the remaining cohomology groups:
\begin{displaymath}
H^k(L) \cong \left\{ \begin{array}{ll}
\Z                                                 & k = 0,2,5,7\\
ker(\pi^*)                                    & k = 3\\
\Z_r                                             & k = 4\\
0                                                  & \textrm{else}.
\end{array} \right.
\end{displaymath}
Observe that $H^3(L) \cong ker(\pi^*)$ will be trivial if and only if $|det( \pi^*)| \neq 0$.

To find $r = |det( \pi^*)|$, we follow the example of \cite{GZ} Proposition 3.3.  Consider the diagram:
\begin{equation} \label{eq:pidiagram} \xymatrix{
	H^3(G) \cong \Z \oplus \Z   &&  H^3(G/K_-^\circ) \oplus H^3(G/K_+^\circ) \cong \Z \oplus \Z \ar[ll]_{ \tau^* \,= \, \tau_-^* - \tau_+^* \qquad \qquad} \\
	H^3(G/H) \ar[u]^{ \eta^*} \cong \Z \oplus \Z   &&  H^3(G/K_-) \oplus H^3(G/K_+) \cong \Z \oplus \Z   \ar[ll]_{ \pi^*\,= \, \pi_-^* - \pi_+^* \qquad \qquad} \ar[u]_{ \mu^*\, = \, \mu_-^* \times \mu_+^*}}
\end{equation}
where the homomorphisms $ \tau_\pm^*$ and $ \eta^*$ are induced by orbit maps, and $ \mu_\pm^*$ are the homomorphisms induced by the maps $gK_\pm^\circ \mapsto gK_\pm$ (which are themselves induced by the inclusions of the identity components $K_\pm^\circ$ in $K_\pm$).  In this case, the orbit $G/K_- = S^3 \times S^2$ is connected, so $ \mu_-^*$ is the identity.  And since $ \mu_+^*$ is an isomorphism by \cite{GWZ} Lemma 13.3a, we have $|det( \mu^*)| = 1$.

We next wish to find $|det( \eta^*)|$.  As previously noted,  $G/H$ is homeomorphic to $S^3 \times (S^3/\langle i \rangle)$.  By uniqueness of the universal cover, the composition $S^3 \times S^3 \xrightarrow{\eta} G/H \xrightarrow{\approx} S^3 \times (S^3/\langle i \rangle)$ induces a homeomorphism of $S^3 \times S^3$ such that the cohomology square:
\begin{equation} \label{eq:etadiagram} \xymatrix{
H^3(S^3 \times S^3)  \cong \Z \oplus \Z & H^3(S^3 \times S^3) \cong \Z \oplus \Z  \ar[l]_{\cong}\\
H^3(G/H) \cong \Z \oplus \Z \ar[u]^{ \eta^*} & H^3(S^3 \times S^3/\langle i \rangle) \cong \Z \oplus \Z \ar[u]_{(id_{S^3} \times f)^*}  \ar[l]_{\cong}}
\end{equation}
commutes, where $f$ is the projection of universal cover of $S^3/\langle i \rangle$ by $S^3$.  An argument involving the K\"unneth isomorphism shows that there are bases for $H^3(S^3 \times S^3/\langle i \rangle)$ and $H^3(S^3 \times S^3)$ such that $(id_{S^3} \times f)^* = {id_{S^3}}^* \times f^*$.  The covering degree $deg(f) = \pm 4$ implies that $|det( \eta^*)| = |det({id_{S^3}}^* \times f^*)| = 4$.

The determinant of $ \tau^*$ follows from \cite{GZ} Proposition 3.3.  They find a basis of $H^3(S^3 \times S^3)$ with respect to which $im( \tau_\pm^*) = \langle(-q_\pm^2, p_\pm^2)\rangle$.  Hence, the absolute value of the determinant of $\tau^* =  \tau_-^* - \tau_+^*$ is $|p_+^2q_-^2 - p_-^2q_+^2|$.

We conclude $|det( \pi^*)| = |det({ \eta^*}^{-1})||det( \tau^*)||det( \mu^*)| = \frac{1}{4}|p_+^2q_-^2 - p_-^2q_+^2|$.  The requirements that $p_+$ and $q_+$ be odd and $p_-, q_- \equiv 1  \textrm{ mod }4$ force $\frac{1}{4}|p_+^2q_-^2 - p_-^2q_+^2|$ to be even, so $H^4(L)$ is a non-trivial cyclic group of even order, and is finite so long as $|p_+^2q_-^2 - p_-^2q_+^2| \neq 0$.  Although the even order of $H^4(L)$ prevents $L$ from being an Eschenburg space, we will show in Section \ref{section4} that it has cohomology type $E_r$ for $r = \frac{1}{4}|p_+^2q_-^2 - p_-^2q_+^2|$.

\smallskip
\noindent \textbf{Case 2.}  On the other hand, suppose $p_+$ is even.  Let $K' $ be the subgroup $\{(e^{jp_+\theta},e^{jq_+\theta})\} \cdot \langle(1,-1), (i,i)\rangle \leq S^3 \times S^3$.  The inclusion of $K_+ = \{(e^{jp_+\theta},e^{jq_+\theta})\} \cdot \langle(i,i)\rangle $ in $K'$ as a subgroup induces a continuous bijection (since $p_+$ is even) from the compact space $G/K_+$ to the Hausdorff space $G/K'$.  It follows that $G/K_+$ is homeomorphic to $G/K'$.  Thus, the cohomology of $G/K_+$ is the same as that of $G/K'$, which was shown in \cite{GWZ} (Lemma 13.6b) to be:
\begin{displaymath}
H^k(G/K_+) \cong \left\{ \begin{array}{ll}
\Z                           & k = 0\\
\Z_4                      & k = 2\\
\Z \oplus \Z_2      & k = 3\\
\Z_2                      & k = 5\\
0                            & \textrm{else}.
\end{array} \right.
\end{displaymath}
Using Sequence \eqref{eq:mv}, Poincar\'e duality and the universal coefficient theorem, we find that the non-trivial cohomology groups of a member $L$ of this subfamily are:  $H^0(L) = H^2(L) = H^7(L) \cong \Z$, $H^5(L) = \Z \oplus \Z_2$, and $H^3(L)$ and $H^4(L)$ are, respectively, the kernel and cokernel of the homomorphism $\pi^* = \pi_-^* - \pi_+^*$ from $H^3(G/K_-) \oplus H^3(G/K_+)$ to $H^3(G/H)$ in Sequence \eqref{eq:mv}.  By Lemma \ref{cycliclemma} (with $t = 2$ and $\kappa = 4$), $H^4(L)$ is cyclic, with order $r = |det(\pi^*|_{\Z \oplus \Z})|$.  In this case, there is a diagram:
\begin{displaymath} \xymatrix{
	H^3(G) \cong \Z \oplus \Z   &&  H^3(G/K_-^\circ) \oplus H^3(G/K_+^\circ) \cong \Z \oplus \Z \ar[ll]_{ \tau^* \,= \, \tau_-^* - \tau_+^* \qquad \qquad} \\
	H^3(G/H) \ar[u]^{ \eta^*} \cong \Z \oplus \Z   &&  H^3(G/K_-) \oplus H^3(G/K_+) \cong \Z \oplus (\Z \oplus \Z_2).  \ar[ll]_{ \pi^*\,= \, \pi_-^* - \pi_+^* \qquad \qquad} \ar[u]_{ \mu^*\, = \, \mu_-^* \times \mu_+^*}}
\end{displaymath}
Comparing this to Diagram \eqref{eq:pidiagram}), we see that the homomorphisms $\eta^*$, $\tau^*$ and $\mu_-^*$ are the same.  By \cite{GWZ} Lemma 13.6, the homomorphism $\mu_+^*$ is multiplication by $\pm4$ on the free part of $H^3(G/K_+)$, while the $\Z_2$ summand is clearly in the kernel.  We conclude that $r = |p_+^2q_-^2 - p_-^2q_+^2|$.  Since $p_+$ is even while $q_\pm$ and $p_-$ are odd, $r$ is always odd, so $H^4(L)$ is finite.  Also, $r = |(p_+q_- + p_-q_+)(p_+q_- - p_-q_+)| \neq 1$ since the parameters $p_\pm$ and $q_\pm$ are non-zero.  Hence, $H^4(L)$ is a non-trivial finite cyclic group of odd order; and by Poincar\'e duality and the universal coefficient theorem, $H^3(L) \cong \Z_2$.  Thus, the cohomology groups of $L$ are:
\begin{displaymath}
H^k(L) \cong \left\{ \begin{array}{ll}
\Z                                   & k = 0,2,7\\
\Z _2                             & k = 3\\
\Z_r                               & k = 4\\
\Z \oplus \Z _2             & k = 5\\
0                                    & \textrm{else}.
\end{array} \right.
\end{displaymath}
Observe that a member of the family $L_{(p_-,q_-),(p_+,q_+)}$ with $p_+$ even shares with members of the family $M_{(p_-,q_-),(p_+,q_+)}$ the distinction of not having cohomology type $E_r$.  The cohomology groups of members of the family $L_{(p_-,q_-),(p_+,q_+)}$ with $p_+$ even differ from those of a type $E_r$ manifold by the presence of $\Z_2$-summands in the third and fifth groups.

\subsection{Members of the family $M_{(p_-,q_-),(p_+,q_+)}$.}

The topology of the family $M_{(p_-,q_-),(p_+,q_+)}$ is described in \cite{GWZ} Theorem 13.1, where they are shown to be 2-connected.  Therefore, they have trivial second and fifth cohomology groups, and cannot be of cohomology type $E_r$.  As a consequence, the Kreck-Stolz invariants cannot be applied to determine diffeomorphism types of these manifolds.

\subsection{Members of the family $N_{(p_-,q_-),(p_+,q_+)}$.}

The cohomology groups of a member $N$ of this family were calculated in \cite{GWZ} Theorem 13.5, and shown to be consistent with those of a manifold having cohomology type $E_r$ where the order of the cyclic group $H^4(N)$ is $r = |p_-^2q_+^2 - p_+^2q_-^2|$.  Since $p_+$ is required to be even while $p_-, q_-$ and $q_+$ are odd, $r$ must be odd and (as in the case of the family $L_{(p_-,q_-),(p_+,q_+)}$ for $p_+$ even) cannot equal one.  Thus, $H^4(N)$ is a non-trivial finite cyclic group of odd order, and the cohomology groups of all members of this family agree with those of Eschenburg spaces.

\subsection{Members of the family $O_{(p,q:m)}$.}

Recall that this family is described by the groups:
\begin{displaymath}
H = \Z_m \,\subseteq \,K_- = \{(e^{ip\theta},e^{iq\theta})\}, \,K_+ = \Delta S^3 \cdot H \,\subseteq \,G = S^3 \times S^3
\end{displaymath}
where $\Delta S^3$ is the diagonal embedding.  The integers $p$ and $q$ are relatively prime, and either $m=1$ (in which case $H$ is the trivial group, and there are no restrictions on the parameters), or $m=2$ (in which case $H = \langle(1,-1)\rangle$ is isomorphic to $\Z_2$ and $p$ is required to be even).  This family naturally splits into two subfamilies, depending on the value of $m$.  In both cases, the non-principal orbit $G/K_- $ is homeomorphic to $S^3 \times S^2$; the difference lies in the other non-principal orbit $G/K_+$, and the principal orbit $G/H$.

\smallskip
\noindent \textbf{Case 1.}  First, suppose $m = 1$.  Then $G/K_+ = S^3 \times S^3/\Delta S^3$ is homeomorphic to $S^3$ under the map sending the coset $[(q_1,q_2)]$ to $q_1{q_2}^{-1}$.  The principal orbit is $G/H = S^3 \times S^3$.  For a member $O_1$ of this subfamily, recall that $H^0(O_1) \cong H^7(O_1) \cong \Z$.  Using Sequence \eqref{eq:mv} and Lemma \ref{cycliclemma} (with $t = 2$ and $\kappa = 4$), one easily sees that the cohomology groups are:
\begin{displaymath}
H^k(O_1) \cong \left\{ \begin{array}{ll}
\Z                                          & k = 0,2,5,7\\
ker( \pi^*)                            & k = 3\\
\Z_r                                      & k = 4\\
0                                           & \textrm{else}.
\end{array} \right.
\end{displaymath}
where again $r = |det( \pi^*)|$ for $ \pi^*$ the homomorphism from the rank two free abelian group $H^3(G/K_-) \oplus H^3(G/K_+)$ to the rank two free abelian group $H^3(G/H)$.  Observe that, in order for $H^3(O_1)$ to be trivial, the determinant of $ \pi^*$ must be non-zero.

As a preliminary step to finding this determinant, let $v$ be a generator of $H^3(S^3)$.  Fix a basis $u_1, u_2$ of $H^3(G/H) = H^3(S^3 \times S^3)$ which corresponds to the images of $v \otimes \1 $ and $\1 \otimes v$ under the K\"unneth isomorphism; that is, $u_i = p_i^*(v)$ where $p_i$ is the projection of the $i^{th}$ factor of $S^3 \times S^3$ onto $S^3$ ($i = 1,2$).  Up to sign, this is the basis used in \cite{GZ} Proposition 3.3 to show that $im( \pi_-^*) = \langle(-q^2, p^2)\rangle$.  We now find $im(\pi_+^*) \leq H^3(S^3 \times S^3)$ with respect to the basis $u_1,u_2$.

Let $S^3 \xhookrightarrow{\Delta} S^3 \times S^3 \xrightarrow{\pi_+} G/K_+ \approx (S^3 \times S^3)/\Delta S^3 \approx S^3$ be the principal $S^3$-bundle with fiber inclusion $\Delta$ the diagonal embedding of $S^3$ in $S^3 \times S^3$.  The composition $\pi_+ \circ \Delta$ is constant, and so is a degree zero map; it follows that the induced homomorphism $\Delta^* \circ \pi_+^*$ is the trivial homomorphism from $H^3(S^3)$ to itself.  Therefore, the image of $\pi_+^*$ is contained in the kernel of $\Delta^*$.  If $\sigma \in C_3(S^3)$ is a singular 3-chain, then for $i = 1,2$:
\begin{displaymath}
\Delta^*(u_i)(\sigma) = u_i(\Delta(\sigma)) = p_i^*(v)((\sigma,\sigma)) = v(\sigma).
\end{displaymath}
So the kernel of $\Delta^*$ is the subgroup of $H^3(S^3 \times S^3)$ generated by $u_1 - u_2$, and there is an integer $n$ such that $im(\pi_+^*) = \langle n(u_1 - u_2) \rangle$.

Next, consider the Serre spectral sequence $(E,d)$ of the Borel fibration $S^3 \times S^3 \xrightarrow{\pi_+} G/K_+ \xrightarrow{\rho} \H P^\infty$ (here, $\rho$ is the classifying map of the previous $S^3$-bundle).  The differential $E_4^{0,3} \cong H^3(S^3 \times S^3) \xrightarrow{d_4} E_4^{4,0} \cong H^4(\H P^{\infty})$ is identified with the transgression (see \cite{McCleary}, Theorem 6.83).  By the definition of the transgression (see, for example, \cite{McCleary} p.186), in this particular instance we have $ker(d_4) = im(\pi_+^*) = \langle n(u_1 - u_2) \rangle$.  Based on the convergence of the spectral sequence to $H^*(G/K_+) \cong H^*(S^3)$, it not difficult to see that $H^3(S^3 \times S^3)/ker(d_4)$ is isomorphic to $H^4(\H P^{\infty}) \cong \Z$.  Using the basis $u_1, u_1 - u_2$ for $H^3(S^3 \times S^3) \cong \Z \oplus \Z$, we conclude that $|n| = 1$; so the image of $\pi_+^*$ in $H^3(S^3 \times S^3)$ with respect to the basis $u_1, u_2$ is the subgroup $\langle (1,-1)\rangle$.

From the above, the absolute value of the determinant of $ \pi^* \, = \, \pi_-^* - \pi_+^*$ is $|p^2 - q^2|$.  The only way that $det(\pi^*)$ can equal zero is if $|p| = |q| = 1$.  Excluding those cases, $H^3(O_1)$ is trivial and $H^4(O_1)$ is a non-trivial finite cyclic group (non-trivial, since $p,q\neq0$ means $|p^2-q^2| = |(p+q)(p-q)| \neq 1$).  So for $|p|$ and $|q|$ not both equal to 1, the cohomology groups of $O_1$ agree with those of a manifold having cohomology type $E_r, \, r = |p^2 - q^2|$.  If either $p$ or $q$ is even, then the order of $H^4(O_1)$ is odd and the cohomology groups of $O_1$ are the same as those of an Eschenburg space.

\smallskip
\noindent \textbf{Case 2.}  Let $m = 2$.  In this case, $G/K_+ = (S^3 \times S^3)/(\Delta S^3 \cdot \langle(1,-1)\rangle)$ is homeomorphic to $\R P^3$ under the map sending the coset $[(q_1,q_2)]$ to the coset $[q_1{q_2}^{-1}]$.  The principal orbit is $G/H = S^3 \times S^3/\langle(1,-1)\rangle$ is homeomorphic to $S^3 \times \R P^3$ under the map $[q_1,q_2] \mapsto (q_1,[q_2])$.  Once again, Sequence \eqref{eq:mv} and Lemma \ref{cycliclemma} (with $t = 2$ and $\kappa = 4$) are sufficient tools for determining the cohomology groups of a member $O_2$ of this subfamily.  As in the previous case, they are:
\begin{displaymath}
H^k(O_2) \cong \left\{ \begin{array}{ll}
\Z                                          & k = 0,2,5,7\\
ker( \pi^*)                            & k = 3\\
\Z_r     & k = 4\\
0                                          & \textrm{else}.
\end{array} \right.
\end{displaymath}
for $r$ the absolute value of the determinant of the homomorphism $\pi^*$ from the rank two free abelian group $H^3(G/K_-) \oplus H^3(G/K_+)$ to the rank two free abelian group $H^3(G/H)$, and $H^3(O_2) \cong ker( \pi^*)$ is trivial when the determinant of $\pi^*$ is not zero.

To calculate $|det( \pi^*)|$, we refer again to Diagram \eqref{eq:pidiagram}.  As before, $ \mu_-$ is the identity map.  Now, however, $ \mu_+$ is the projection of the universal cover of $\R P^3$ by $S^3$, which has covering degree two; so $|det( \mu^*)| = 2$.  The composition  $S^3 \times S^3 \xrightarrow{\eta} G/H \xrightarrow{\approx} S^3 \times \R P^3$ is the universal cover, and an argument analogous to that involving Diagram \eqref{eq:etadiagram} shows that $|\det( \eta^*)| = 2$.

The absolute value of the determinant of $ \tau^*$ has already been computed; the homomorphism $ \pi^*$ that determined the order of the fourth cohomology group in the previous subfamily $O_1$ is the same as the current homomorphism $ \tau^*$.  Thus,  the absolute value of the determinant of the current homomorphism $ \pi^*$ is $|det(\pi^*)| = |det({ \eta^*)}^{-1}||det( \tau^*)||det( \mu^*)| = |p^2 - q^2|$.  Recall that in this case $p$ is even, so $H^4(O_2)$ is finite cyclic of odd order $r = |p^2 - q^2|$, $H^3(O_2) \cong ker( \pi^*) = 0$, and the cohomology groups of $O_2$ are the same as those of Eschenburg spaces.

\section{Proof of Proposition 1.1, continued:  the cohomology rings.} \label{section4}

We have shown that the cohomology groups of members of the family $L_{(p_-,q_-),(p_+,q_+)}$ where $p_+$ is odd and $p_+^2q_-^2 - p_-^2q_+^2 \neq 0$, all members of the family $N_{(p_-,q_-),(p_+,q_+)}$, and members of the family $O_{(p,q:m)}$ for $|p|$ and $|q|$ not both equal to one, are in agreement with those of a manifold having cohomology type $E_r$.  Now, we show that if $M$ is any of the above manifolds, then the cohomology ring $H^*(M)$ is generated by classes $x \in H^2(M)$ and $y \in H^5(M)$.  Not only does this show that all of the above manifolds have cohomology type $E_r$; it also implies that their cohomology rings are the same as those of Eschenburg spaces whenever the order $r$ of the fourth cohomology group is odd.

We also give an almost complete description of the cohomology ring structure of the remaining manifolds, the members of $L_{(p_-,q_-),(p_+,q_+)}$ such that $p_+$ is even.  We show that, for any such manifold $M$ and $x$ and $y$ generators of $H^2(M)$ and the free part of $H^5(M)$ respectively, the class $x^2$ generates $H^4(M)$ and $xy$ generates $H^7(M)$.

For all of the above, the non-principal orbits $G/K_-$ are closed orientable submanifolds of codimension two.  The manifolds themselves are orientable (since simply connected), so the normal disk bundles over $G/K_-$ are orientable bundles with fiber $D^2$.  Thus, we have at our disposal Sequence \eqref{eq:les} and (provided the conditions are met) Lemma \ref{genlemma}, setting $t = 2$ in both.  In all that follows, we will assume that the class $x$ generates $H^2(M)$, the class $y$ generates (the free part of) the $H^5(M)$, and $\1 _\pm$ is the unit of $H^*(G/K_\pm)$.

\subsection{Cohomology ring of $L_{(p_-,q_-),(p_+,q_+)}$.}

\noindent \textbf{Case 1.}  Let $L$ be a member of the subfamily of $L_{(p_-,q_-),(p_+,q_+)}$ for which $p_+$ is odd and $p_+^2q_-^2 - p_-^2q_+^2 \neq 0$.  In this case, Lemma \ref{genlemma} cannot be called on to show that the square of the generator $x$ of $H^2(L)$ generates $H^4(L)$, as Condition \ref{condition3} fails.  To see why, we analyze the section of Diagram \eqref{eq:ladder} corresponding to $k = 1$ and $2$.

Since $H^1(G/K_-)$ and $H^1(G/H)$ are both trivial, the homomorphisms $j_-^*$  from $H^2(L,G/K_-)$ to $H^2(L) \cong \Z$ and $j^*$ from $H^2(D(G/K_+),G/H)$ to $H^2(G/K_+) \cong \Z_2$ are injective.  Because the groups $H^2(L,G/K_-)$ and $H^2(D(G/K_+),G/H)$ are isomorphic, injectivity of the homomorphisms $j_-^*$ and $j^*$ implies that $H^2(L,G/K_-)$ and $H^2(D(G/K_+),G/H)$ are trivial.  Since $H^3(L) = 0$, the isomorphism of $H^3(L,G/K_-)$ and $H^3(D(G/K_+),G/H)$ (together with commutativity of the diagram) implies that the homomorphism $j^*$ from $H^3(D(G/K_+),G/H)$ to $H^3(G/K_+)$ is the zero homomorphism.  This gives two short exact sequences:
\begin{displaymath}
\begin{array}{c}
0 \rightarrow H^2(L) \cong \Z \cdot x \xrightarrow{i_-^*} H^2(G/K_-) \cong \Z \cdot \gamma \xrightarrow{\delta_-} H^3(L,G/K_-) \rightarrow 0 \\
0 \rightarrow H^2(G/K_+) \cong \Z_2 \xrightarrow{i^*} H^2(G/H) \cong \Z_4 \xrightarrow{\delta} H^3(D(G/K_+),G/H) \rightarrow 0
\end{array}
\end{displaymath}
where the groups $H^3(L,G/K_-)$ and $H^3(D(G/K_+), G/H)$ are isomorphic.  From the second sequence, it follows that $H^3(D(G/K_+),G/H)$ is isomorphic to $\Z_2$.  Then, by the first sequence, the homomorphism $i_-^*$ from $H^2(L)$ to $H^2(G/K_-)$ is multiplication by 2 and $i_-^*(x) = \pm 2\gamma$ for $\gamma$ a generator of $H^2(G/K_-)$.  So if Condition \ref{condition3} of Lemma \ref{genlemma} holds, then the order of $H^4(L)$ must be relatively prime to $2$.  We have already shown, however, that the order of $H^4(L)$ is even.

Fortunately, there is another way to see that $x^2$ generates $H^4(L)$.  Setting $t = 2$, $k = 2$ and $B_\pm = G/K_\pm$ in Sequence \eqref{eq:les} for the pair $(L,G/K_+)$ gives a short exact sequence:
\begin{displaymath}
0 \rightarrow H^0(G/K_-) \cong \Z \xrightarrow{J} H^2(L) \cong \Z \cdot x \xrightarrow{i_+^*} H^2(G/K_+) \cong \Z_2 \rightarrow 0
\end{displaymath}
and we see that $J(\1_-) = \pm 2x$.  Setting $k = 4$ in Sequence \eqref{eq:les}, exactness together with the triviality of $H^4(G/K_+)$ implies $J(\gamma)$ generates $H^4(L)$.  Then:
\begin{displaymath}
2J(\gamma) = J(2\gamma) = J(i_-^*(x)) = J(\1_-) \smile x = \pm 2x^2
\end{displaymath}
(recall that $J$ is an $H^*(L)$-module homomorphism).  Since $J(\gamma)$ generates $H^4(L)$, the subgroup generated by $2x^2 = \pm 2J(\gamma)$ is an index two subgroup, and so is maximal in $H^4(L)$.  We show that $x^2$ is not an element of $\langle 2x^2 \rangle$, from which it follows $x^2$ generates $H^4(L)$.  This argument will require both integral and $\Z_2$ cohomology, so we temporarily resort to explicitly indicating coefficients.

The short exact sequence of abelian groups:
\begin{displaymath}
0 \rightarrow \Z \xrightarrow{h} \Z \xrightarrow{g} \Z_2 \rightarrow 0
\end{displaymath}
where $h$ is multiplication by two and $g$ is the natural projection, gives rise to the long exact cohomology sequence (see \cite{Massey}):
\begin{displaymath}
\cdots \rightarrow H^k(L; \Z) \xrightarrow{h_\#} H^k(L;\Z) \xrightarrow{g_\#} H^k(L;\Z_2) \xrightarrow{\beta} H^{k+1}(L;\Z) \rightarrow \cdots
\end{displaymath}
where $\beta$ is the Bockstein operator.  Since $\pm 2x^2 =  2J(\gamma) = h_\#(J(\gamma))$, exactness implies $\langle 2x^2 \rangle$ is the kernel of $g_\#$.  If $g_\#(x^2)$ can be shown to be non-trivial in $H^4(L;\Z_2)$, it will follow that $x^2$ is not in $\langle 2x^2 \rangle$.

Since $H^3(L;\Z)$ is trivial, exactness implies that the homomorphism $g_\#$ from $H^2(L;\Z)$ to $H^2(L;\Z_2)$ is surjective.  Thus $g_\#$ sends the generator $x$ of $H^2(L;\Z)$ to the generator $w$ of $H^2(L;\Z_2)$.  Checking the definitions of the induced homomorphism $g_\#$ and the cup product reveals that $g_\#(x^2) = w^2$.

The $\Z_2$-cohomology of $L$ and the non-principal orbits $G/K_\pm$ are:
\begin{displaymath}
\begin{array}{c}
H^k(L;\Z_2) \cong \left\{ \begin{array}{ll}
\Z_2                                        & \textrm{if }k = 0,2,3,4,5,7\\
0                                                              & \textrm{else}
\end{array} \right.\\\\

H^k(G/K_-; \Z_2) \cong \left\{ \begin{array}{ll}
\Z_2                                            & \textrm{if }k = 0,2,3,5\\
0                                                                  & \textrm{else}
\end{array} \right.\\\\

H^k(G/K_+; \Z_2) \cong \left\{ \begin{array}{ll}
\Z_2                                  & \textrm{if }k = 0,1,2,3,4,5\\
0                                                        & \textrm{else}.
\end{array} \right.
\end{array}
\end{displaymath}
Applying Sequence \eqref{eq:les} with $\Z_2$-coefficients to the pair $(L,G/K_-)$ (recall that $G/K_+$ need not be orientable if we take coefficients in $\Z_2$) reveals that $J$ is an isomorphism from $H^0(G/K_+; \Z_2)$ to $H^2(L; \Z_2)$.  Let $\1 $ be the unit of the cohomology ring $H^*(G/K_+; \Z_2)$; then $J(\1 ) = w$.  The corresponding long exact sequence for the pair  $(L,G/K_+)$ shows the homomorphism $i_+^*$ from $H^2(L; \Z_2)$ to $H^2(G/K_+; \Z_2)$ is an isomorphism, so $i_+^*(w)$ generates $H^2(G/K_+; \Z_2)$.  Returning to the sequence of the pair $(L,G/K_-)$, we see that $H^2(G/K_+; \Z_2)$ is isomorphic to $H^4(L; \Z_2)$ under $J$; hence, $J(i_+^*(w))$ generates $H^4(L; \Z_2)$.  Since $J$ is an $H^*(L;\Z_2)$-module homomorphism:
\begin{displaymath}
J(i_+^*(w)) = J(\1 \smile i_+^*(w)) = J(\1 ) \smile w = w^2.
\end{displaymath}
Thus, $w^2 = g_\#(x^2)$ generates $H^4(L:\Z_2)$.  In particular, $g_\#(x^2)$ is non-trivial, which is what we needed to show in order to conclude that $x^2$ generates $H^4(L;\Z)$.  This is the last time $\Z_2$ coefficients will be used, and we return to the earlier convention of implicitly assuming integral cohomology.

We now show that all of the conditions of Lemma \ref{genlemma} do hold when $\kappa = 7$, from which it follows that the class $xy$ generates $H^7(L)$.  Recall that $t = 2$, and observe that all conditions on the cohomology groups are met.  Since $H^1(G/K_-)$ is trivial, $i_+^*$ from $H^2(L)$ to $H^2(G/K_+)$ is a surjection.

We check the remaining three conditions.  Since $H^7(G/K_+)$ is trivial, $i_+^*$ from $H^7(L)$ to $H^7(G/K_+)$ is the zero homomorphism, and Condition \ref{condition1} holds.  Condition \ref{condition2} also holds, with $H^5(G/K_-) \cong \Z \cdot \nu$.  The group $H^2(G/K_+)$ is finite cyclic of order two, and $H^7(L)$ is infinite cyclic; so to verify Condition \ref{condition3}, we need to show that $i_+^*(y) = \pm 2 \nu$.

Take $k = 5$ in Diagram \eqref{eq:ladder}.  Since $H^4(G/K_-)$ and $H^4(G/H)$ are trivial, there are injections of $H^5(L,G/K_-)$ into the infinite cyclic group $H^5(L)$ and of $H^5(D(G/K_+),G/H)$ into the finite cyclic group $H^5(G/K_+)$.  Because the groups $H^5(L,G/K_-)$ and $H^5(D(G/K_+),G/H)$ are isomorphic, the only way both injections can hold is if $H^5(L,G/K_-) = H^5(D(G/K_+),G/H) = 0$.  This gives two short exact sequences:
\begin{displaymath}
\begin{array}{c}
0 \rightarrow H^5(L) \cong \Z \cdot y \xrightarrow{i_-^*} H^5(G/K_-) \cong \Z \cdot \nu \xrightarrow{\delta_-} H^5(L,G/K_-) \rightarrow 0 \\
0 \rightarrow H^5(G/K_+) \cong \Z_2 \xrightarrow{i^*} H^5(G/H) \cong \Z_4 \xrightarrow{\delta} H^5(D(G/K_+),G/H) \rightarrow 0
\end{array}
\end{displaymath}
where the groups $H^5(L,G/K_-)$ and $H^5(G/K_+, G/H)$ are isomorphic.  From the second sequence, we conclude $H^6(D(G/K_+),G/H) \cong \Z_2$.  It is then apparent by the first sequence that $i_-^*$ is multiplication by two; so $i_-^*(y) = \pm 2\nu$, Condition \ref{condition3} holds, and by Lemma \ref{genlemma}, $xy$ generates $H^7(L)$.

We conclude that a member of the family $L_{(p_-,q_-),(p_+,q_+)}$ for which $p_+$ is odd and $p_+^2q_-^2 - p_-^2q_+^2 \neq 0$ has cohomology type $E_r, \, r = \frac{1}{4}|p_+^2q_-^2 - p_-^2q_+^2|$, proving Proposition \ref{typeErprop} (a).

\smallskip
\noindent \textbf{Case 2.}  Suppose $L$ is a member of the subfamily of $L_{(p_-,q_-),(p_+,q_+)}$ such that $p_+$ is even.  Let $x$ generate $H^2(L) \cong \Z$ and $y$ the free part of $H^5(L) \cong \Z \oplus \Z_2$; we show that $x^2$ generates $H^4(L)$ and $xy$ generates $H^7(L)$.  Whether or not the classes $x$ and $y$, together with a class $\xi$ generating $H^3(L) \cong \Z_2$, form a complete set of generators for the ring $H^*(L)$ is unknown at this time.

Setting $t = 2$, we confirm that the conditions of Lemma \ref{genlemma} hold when $\kappa = 4,7$.  The conditions on $H^2(L)$ and $H^2(G/K_+)$, as well as Conditions \ref{condition1} and \ref{condition2}, are easily checked.  Sequence \eqref{eq:les} can be used to verify that the inclusion-induced homomorphism from $H^2(L) \cong \Z$ to $H^2(G/K_+) \cong \Z_4$ is a surjection, and those from $H^\kappa(L)$ to $H^\kappa(G/K_+)$ (for $\kappa = 4,7$) are multiplication by zero.  It remains to check Condition \ref{condition3}.

When $\kappa = 4$, $H^4(L)$ is finite cyclic.  We show $i_-^*(x)$ generates $H^2(G/K_-)$ when $x$ generates $H^2(L)$.  Consider Diagram \eqref{eq:ladder}.  The second relative cohomology groups, which are isomorphic, inject into both a free group and a finite group and so must be trivial.  Because $H^3(L) \cong \Z_2$ is in the kernel of $i_-^*$, the third relative group $H^3(L,G/K_-)$ surjects onto $H^3(L)$.  By exactness of the top row, $H^3(L,G/K_-)/im(\delta)$ is isomorphic to $\Z_2$ and $im(\delta)$ is finite cyclic.  We conclude that $H^3(L,G/K_-)$ is a non-trivial finite group.

Triviality of $H^2(D(G/K_+),G/H)$ implies $H^2(G/K_+) \cong \Z_4$ injects into $H^2(G/H) \cong \Z_4$; so $H^2(G/K_+)$ and $H^2(G/H)$ are isomorphic.  It follows from the isomorphism of the third relative cohomology groups, together with exactness in the bottom row, that $H^3(L,G/K_-)$ injects into $H^3(G/K_+) \cong \Z \oplus \Z_2$.  So $H^3(L,G/K_-)$ is isomorphic to $\Z_2$, and the surjection of $H^3(L,G/K_-)$ onto $H^3(L)$ is an isomorphism.  Then, by exactness of the top row, the inclusion-induced homomorphism $i_-^*$ from $H^2(L) \cong \Z\cdot x$ to $H^2(G/K_-)$ must also be an isomorphism, and $i_-^*(x)$ generates $H^2(G/K_-)$.  This satisfies Condition \ref{condition3} in the case $\kappa = 4$.

If $\kappa = 7$, $H^7(L)$ is infinite cyclic.  We show that, for $y$ a generator of the free part of $H^5(L)$, $i_-^*(y)$ is four times a generator of $H^5(G/K_-) \cong \Z$.  Again turning to Diagram \eqref{eq:ladder}, we see that the fifth relative cohomology groups, which are isomorphic, inject into both $H^5(L) \cong \Z \cdot y \oplus \Z_2$ and $H^5(G/K_+) \cong \Z_2$.  Hence, they are either trivial or cyclic of order two.  Since the $\Z_2$ summand of $H^5(L)$ is in the kernel of the homomorphism $i_-^*$ from $H^5(L)$ to $H^5(G/K_-)$, we conclude that the fifth relative cohomology groups are isomorphic to $\Z_2$.  Then exactness of the bottom row together with the isomorphism of the relative groups gives an isomorphism between $H^6(L,G/K_-)$ and $H^5(G/H) \cong \Z_4$.  Restricting $i_-^*$ to the free part of $H^5(L)$ gives rise to a short exact sequence:
\begin{displaymath}
0 \rightarrow \Z \cdot y \xrightarrow{{i_-^*}|_\Z} H^5(G/K_-) \cong \Z \xrightarrow{\delta} \Z_4 \rightarrow 0.
\end{displaymath}
Hence, $i_-^*(y)$ is four times a generator of $H^5(G/K_-)$, satisfying Condition \ref{condition3} in the case $\kappa = 7$, and by Lemma \ref{genlemma} it follows that $x^2$ generates $H^4(L)$ and $xy$ generates $H^7(L)$.

\subsection{Cohomology ring of $N_{(p_-,q_-),(p_+,q_+)}$.}

Let $N$ be a member of the family $N_{(p_-,q_-),(p_+,q_+)}$.  The cohomology groups of the orbits (as computed in \cite{GWZ} Lemma 13.6) are:
\begin{displaymath}
\begin{array}{cc}
H^k(G/K_-) \cong \left\{ \begin{array}{ll}
\Z                                      & k = 0,3,5\\
\Z \oplus \Z_2 & k = 2\\
\Z_2                                  & k = 4\\
0                                                        & \textrm{else}
\end{array} \right.
&
H^k(G/K_+) \cong \left\{ \begin{array}{ll}
\Z                 & k = 0\\
\Z_4            & k = 2\\
\Z \oplus \Z_2 & k = 3\\
\Z_2            & k = 5\\
0                                  & \textrm{else}
\end{array} \right.
\end{array}
\end{displaymath}

\begin{displaymath}
H^k(G/H) \cong \left\{ \begin{array}{ll}
\Z                                                 & k = 0,6\\
\Z_2 \oplus \Z_4       & k = 2,5\\
\Z \oplus \Z \oplus \Z_2 & k = 3\\
\Z_2                                                                  & k = 4\\
0                                                                  & \textrm{else}
\end{array} \right.
\end{displaymath}

The classes $x$ and $y$ respectively generate the infinite cyclic groups $H^2(N)$ and $H^5(N)$.  To show that $x^2$ generates $H^4(N)$ and $xy$ generates $H^7(N)$, we turn to Lemma \ref{genlemma} (recall that $t = 2$).  For $ \kappa = 4$ and $7$, all of the conditions on the cohomology groups are satisfied, including Condition \ref{condition2}.  In particular, $H^2(G/K_+)$ is finite cyclic of order $n=4$.  Taking $k = 2$ in Sequence \eqref{eq:les}, one sees that the inclusion-induced homomorphism from $H^2(N)$ to $H^2(G/K_+)$ is a surjection, and also that the inclusion-induced homomorphisms from $H^\kappa(N)$ to $H^\kappa(G/K_+)$, $\kappa = 4,7$, are the zero homomorphisms (Condition \ref{condition1}).  It remains only to check that the requirements of Condition \ref{condition3} are satisfied.

For $\kappa = 4$, the group $H^4(N)$ is finite cyclic.  Suppose the image of $x$ under the inclusion-induced homomorphism $i_-^*$ is the element $(s,\beta)$ in $H^2(G/K_-) \cong \Z \oplus \Z_2$.  In order for Condition \ref{condition3} to hold, $s$ must be relatively prime to the order of $H^4(N)$.  We claim this is true; that, in fact, $|s| = 1$.  To see this, consider Diagram \eqref{eq:ladder}.  Setting $k = 2$, we see that the second relative cohomology groups, which are isomorphic, inject into both the infinite cyclic group $H^2(N)$ and the finite cyclic group $H^2(G/K_+)$; so they must be trivial.  Because $H^3(N)$ is trivial, the homomorphism $\delta_-$ from $H^2(G/K_-)$ to $H^3(N,G/K_-)$ is a surjection.  Commutativity of the diagram together with the isomorphism of the third relative groups forces the homomorphism $\delta$ from $H^2(G/H)$ to $H^3(D(G/K_+),G/H)$ to be surjective as well.  This gives two short exact sequences:
\begin{displaymath}
\begin{array}{c}
0 \rightarrow H^2(N) \cong \Z \cdot x \xrightarrow{i_-^*} H^2(G/K_-) \cong \Z \oplus \Z_2 \xrightarrow{\delta_-} H^3(N,G/K_-) \rightarrow 0\\
0 \rightarrow H^2(G/K_+) \cong \Z_4 \stackrel{i^*}{\rightarrow} H^2(G/H) \cong \Z_2 \oplus \Z_4 \stackrel{\delta}{\rightarrow} H^3(D(G/K_+),G/H) \rightarrow 0
\end{array}
\end{displaymath}
where the relative cohomology groups are isomorphic.  From the second sequence we see that the order of the relative groups is the order of $H^2(G/H)$ divided by the order of $H^2(G/K_+)$; hence, the relative groups are isomorphic to $\Z_2$.

Now consider the first sequence.  Because $i_-^*$ is injective and $i_-^*(x) = (s,\beta)$, $s$ cannot be zero.  By exactness, $H^3(N,G/K_-) \cong \Z_2$ is isomorphic to $H^2(G/K_-)/im(i_-^*)$.  If $\beta = [0]$, the group $H^2(G/K_-)/ \langle (s,\beta)\rangle$ is clearly isomorphic to $\Z_s \oplus \Z_2$.  If instead $\beta = [1]$, the surjective homomorphism from $H^2(G/K_-) \cong \Z \oplus \Z_2$ to $\Z_{2s}$, defined by $(1,[0]) \mapsto [1]$ and $(0,[1]) \mapsto [s]$, has kernel $\langle (s,[1]) \rangle$; hence, $H^2(G/K_-)/\langle (s,[1])\rangle$ is isomorphic to $\Z_{2s}$.  In both cases, $H^2(G/K_-)/im(i_-^*)$ is a finite group with $2|s|$ elements.  Because we know $H^2(G/K_-)/im(i_-^*)$ is isomorphic to $\Z_2$, we conclude $|s|=1$.  Thus, $s$ is relatively prime to the order of $H^4(N)$, Condition \ref{condition3} is satisfied and Lemma \ref{genlemma} holds for $ \kappa = 4$.  We have shown $x^2$ generates $H^4(N)$.

We now show that Condition \ref{condition3} of Lemma \ref{genlemma} holds for $\kappa = 7$, from which it follows $xy$ generates $H^7(N)$.  Because $H^7(N)$ is infinite cyclic and $H^2(G/K_+)$ is finite cyclic of order $n=4$, Condition \ref{condition3} requires the image of the generator $y$ of $H^5(N)$ under $i_-^*$ to be (up to sign) four times a generator of the infinite cyclic group $H^5(G/K_-)$.

To show this is true, set $k = 5$ in Diagram \eqref{eq:ladder}.  The sixth relative cohomology groups are isomorphic, and by exactness of the bottom row are isomorphic to the quotient of $H^5(G/H)$ by $i^*(H^5(G/K_+))$.  So the orders of the sixth relative groups are equal to the order of $H^5(G/H) \cong \Z_2 \oplus \Z_4$ divided by the order of $i^*(H^5(G/K_+))$.  Since $H^5(G/K_+) \cong \Z_2$, this is either four or eight.  Observe that these relative groups cannot contain elements of order eight, since they are isomorphic to a quotient of $\Z_2 \oplus \Z_4$; therefore, if they are eight element groups, they cannot be cyclic.  But the infinite cyclic group $H^5(G/K_-)$ surjects onto the relative cohomology groups, so we conclude they are cyclic of order four.  It then follows from exactness of the top row that the homomorphism $i_-^*$ from $H^5(N)$ to $H^5(G/K_-)$ is multiplication by four.  Hence, by Lemma \ref{genlemma}, $xy$ generates $H^7(N)$.

We have now shown that any member $N$ of the family $N_{(p_-,q_-),(p_+,q_+)}$ not only has cohomology type $E_r$ where $r = |p_+^2q_-^2 - p_-^2q_+^2|$, but also that $H^*(N)$ has the same cohomology ring as an Eschenburg space.  This proves Propositions \ref{typeErprop} (b) and \ref{eschencor} (a).

\subsection{Cohomology ring of $O_{(p,q:m)}$.}

Let $O_m$ be a member of this family; $H^*(O_m)$ is the most straightforward of the rings to calculate.  Although two cases (arising from differences in the orbits $G/K_+$ and $G/H$ depending on whether $m=1$ or $m=2$) need to be considered separately, calculations are greatly simplified by the orientability of both non-principal orbits.  This means Sequence \eqref{eq:les} holds for both of the pairs $(O_m,G/K_\pm)$.  Despite differences in some of the cohomology groups, arguments for each of the cases $m=1,2$ are similar; we sketch the general method.

For $\kappa = 4,7$, all conditions of Lemma \ref{genlemma} applying to the cohomology groups (including Condition \ref{condition2}) are met.  As before, $t = 2$.  Sequence \eqref{eq:les} (applied to the pair $(O_m,G/K_+)$) can be used to show that $H^2(O_m)$ surjects onto $H^2(G/K_+)$.  This same sequence can be used to show that the homomorphisms $i_+^*$ from $H^{\kappa}(O_m)$ to $H^{\kappa}(G/K_+)$, $\kappa = 4,7$, are trivial homomorphisms, and consequently that Condition \ref{condition1} holds.  Finally, applying Sequence \eqref{eq:les} to the pair $(O_m,G/K_-)$, one sees that under the homomorphism $i_-^*$, the image of a generator of $H^{\kappa - 2}(O_m)$ meets the requirements of Condition \ref{condition3}.  Thus, by Lemma \ref{genlemma}, $O_m$ has cohomology type $E_r, \, r = |p^2 - q^2|$ so long as $|p|$ and $|q|$ are not both equal to one.  This proves Proposition \ref{typeErprop} (c).  If, in addition, either $p$ or $q$ is even, then $r$ is odd, proving Corollary \ref{eschencor} (b).

This completes the proof of Proposition \ref{typeErprop} and Corollary \ref{eschencor}.  \qed

\subsection{Cohomology rings:  summary.}
We conclude with a summary of the currently known integral cohomology rings:\\

\begin{tabular}{l|l}
\hline \hline
\T \B Family & $L_{(p_-,q_-),(p_+,q_+)}$:  $p_+$ odd and  $\frac{p_-}{q_-} \neq \frac{p_+}{q_+}$\\
\hline
\T \B Non-trivial & $H^0(L) = H^2(L) = H^5(L) = H^7(L) = \Z$\\
\B cohomology groups & $H^4(L) = \Z_r $ for $r = \frac{1}{4}|p_+^2q_-^2 - p_-^2q_+^2|$\\
\hline
\T \B Ring generators & $x \in H^2(L)$ and $y \in H^5(L)$\\
\hline
\T \B Notes & Has cohomology type $E_r$, $r$ even.\\
\hline \hline
\end{tabular}

\begin{tabular}{l|l}
\hline \hline
\T \B Family & $L_{(p_-,q_-),(p_+,q_+)}$:  $p_+$ even\\
\hline
\T \B Non-trivial & $H^0(L) = H^2(L) = H^7(L) = \Z$\\
\B cohomology groups & $H^3(L) = \Z_2$,  $H^5(L) = \Z \oplus \Z_2$\\
\B & $H^4(L) = \Z_r $ for $r = |p_+^2q_-^2 - p_-^2q_+^2|$\\
\hline
\T \B Ring generators & Let $H^2(L) = \Z \cdot x$ and $H^5(L) = \Z \cdot y \oplus \Z_2$.\\
\B (partial list) & Then $x^2$ generates $H^4(L)$\\
\B & and $xy$ generates $H^7(L)$.\\
\hline
\T \B Notes & $r$ is always odd.\\
\hline \hline
\end{tabular}

\begin{tabular}{l|l}
\hline \hline
\T \B Family & $M_{(p_-,q_-),(p_+,q_+)}$:  $\frac{p_-}{q_-} \neq \frac{p_+}{q_+}$\\
\hline
\T \B Non-trivial& $H^0(M) = H^7(M) = \Z$\\
\B cohomology groups & $H^4(M) = \Z_r$ for $r = \frac{1}{8}|p_+^2q_-^2 - p_-^2q_+^2|$\\
\hline
\T \B Ring generators & $y \in H^4(M)$ and $z \in H^7(M)$\\
\hline
\T \B Notes & Computed in \cite{GWZ}.\\
    & Has the same cohomology ring as an\\
\B & $S^3$-bundle over $S^4$.\\
\hline \hline
\end{tabular}

\begin{tabular}{l|l}
\hline \hline
\T \B Family & $N_{(p_-,q_-),(p_+,q_+)}$\\
\hline
\T \B Non-trivial& $H^0(N) = H^2(N) = H^5(N) = H^7(N) = \Z$\\
\B cohomology groups & $H^4(N) = \Z_r$ for $r = |p_+^2q_-^2 - p_-^2q_+^2|$\\
\hline
\T \B Ring generators & $x \in H^2(N)$ and $y \in H^5(N)$\\
\hline
\T \B Notes & Groups computed in \cite{GWZ}.\\
\B & Has cohomology type $E_r$, $r$ odd.\\
\hline \hline
\end{tabular}

\begin{tabular}{l|l}
\hline \hline
\T \B Family & $O_{(p,q:m)}$:  either $|p|$ or $|q|$ to not equal to 1.\\
\hline
\T \B Non-trivial& $H^0(O) = H^2(O) = H^5(O) = H^7(O) = \Z$\\
\B cohomology groups & $H^4(O) = \Z_r$ for $r = |p^2 - q^2|$\\
\hline
\T \B Ring generators & $x \in H^2(O)$ and $y \in H^5(O)$\\
\hline
\T \B Notes & Has cohomology type $E_r$.\\
\B & If either $p$ or $q$ is even, $r$ is odd.\\
\hline \hline
\end{tabular}


\end{document}